\documentclass[letterpaper, reqno,11pt]{article}
\usepackage[margin=1.0in]{geometry}
\usepackage[percent]{overpic}

\usepackage{latexsym,color,graphicx,amsmath,amssymb,amsthm, enumitem, cancel, color, subfigure}

\def\R{{\mathbb R}}

\def\p{{\bf p}}
\def\q{{\bf q}}

\def\rank{{\rm rank}}

\title{Distinct distances for points lying on curves in $\R^d$---the bipartite case 
}

\author{
Hadas Baer-Erenfeld\thanks{%
School of Computer Science and Engineering, The Hebrew University of Jerusalem,
Jerusalem, Israel.
{\sl hadas.baer@mail.huji.ac.il}}~ and Orit E. Raz\thanks{%
Einstein Institute of Mathematics, The Hebrew University of Jerusalem,
Jerusalem, Israel.
{\sl oritraz@mail.huji.ac.il} }
}
\newtheorem{theorem}{Theorem}[section]

\newtheorem{corollary}[theorem]{Corollary}
\newtheorem{conj}[theorem]{Conjecture}
\newtheorem{lemma}[theorem]{Lemma}

\theoremstyle{definition}
\newtheorem{defin}[theorem]{Definition}
\newtheorem{remark}[theorem]{Remark}

\DeclareMathOperator{\im}{Im}

\begin{document}
\maketitle

\begin{abstract}
Let $\gamma_1,\gamma_2$ be a pair of constant-degree irreducible algebraic 
curves in $\R^d$.
Assume that $\gamma_i$ is neither contained in a hyperplane nor in a quadric surface in $\mathbb{R}^d$, for each $i=1,2$.
We show that for every pair of $n$-point sets $P_1\subset\gamma_1$ and $P_2\subset\gamma_2$, 
the number of distinct distances spanned by $P_1\times P_2$ is $\Omega(n^{3/2})$, with 
a constant of proportionality that depends on $\deg\gamma_1$, $\deg\gamma_2$, and $d$. 
This extends earlier results of Charalambides~\cite{Char}, Pach and De Zeeuw~\cite{PdZ}, and Raz~\cite{Razdist} to the bipartite version.
For the proof we use rigidity theory, and in particular the description of Bolker and Roth~\cite{BolkerRoth1980} for realizations in $\R^d$ of the complete bipartite graph $K_{m,n}$ that are not infinitesimally rigid. 
\end{abstract}

\section{Introduction}
The distinct distances problem of Erd\H os~\cite{Erd} asks for the minimum number of distinct 
distances spanned by any set $P$ of $n$ points in the plane. The $\sqrt n\times\sqrt n$ integer 
grid in the plane induces $\Theta (n/\sqrt{\log n})$ distinct distances, and Erd\H os conjectured
that this number is asymptotically tight. In 2010, Guth and Katz~\cite{GK2} proved that every set 
of $n$ points in the plane spans $\Omega(n/\log n)$ distinct distances, which almost matches Erd\H os's upper bound.

An instance of the problem suggested by Purdy (e.g., see \cite[Section 5.5]{BMP})
asks for the minimum number of distinct distances spanned by
pairs of $P_1\times P_2$, where, for $i=1,2$, $P_i$ is a set of $n$ points that lie on a line $\ell_i$ in the plane.
Elekes and R\'onyai \cite{ER00} showed that, in contrast with the general case, this number is $\omega(n)$, unless
the lines $\ell_1, \ell_2$ are either orthogonal or parallel to one another 
(where in the latter cases sets with only $O(n)$ distinct distances between them can easily be constructed).
Sharir, Sheffer, and Solymosi~\cite{SSS} strengthened the result by showing that 
the number of distinct distances spanned by $P_1\times P_2$ in the non-parallel, non-orthogonal case 
is $\Omega(n^{4/3})$. This result was later generalized by Pach and De Zeeuw~\cite{PdZ} to 
the case where, for $i=1,2$, $P_i$ is a set of points that lie on some irreducible constant-degree
algebraic curve $\gamma_i$ in the plane. They showed that in this case the number of distinct 
distances spanned by $P_1\times P_2$ is again $\Omega(n^{4/3})$, unless $\gamma_1,\gamma_2$ are
either a pair of orthogonal lines, a pair of (possibly coinciding) parallel lines, or a pair 
of (possibly coinciding) concentric circles.

The first to consider the distinct distances problem (in the plane and in higher dimensions), 
with points restricted to a single arbitrary constant-degree algebraic curve, was Charalambides~\cite{Char}. 
He showed that if a set $P$ of $n$ points is contained in a constant-degree algebraic curve $\gamma$
in $\R^d$, for any $d\ge 2$, then the number of distinct distances spanned by $P$ is $\Omega(n^{5/4})$, 
unless $\gamma$ is an {\it algebraic helix} (see Charalambides~\cite[Definition~1.5 and Lemma~7.4]{Char}). Raz~\cite{Razdist} 
has shown that the bound $\Omega(n^{5/4})$ in \cite{Char} can be improved to $\Omega(n^{4/3})$,
also for $d>2$. In the plane, an algebraic helix is just a line or a circle, so the result of 
Pach and De Zeeuw provides a generalization (to the bipartite case)
and an improved bound of Charalambides' result for the case $d=2$. 
A recent general result of Solymosi--Zahl~\cite{SolyZahl}, implies in particular that the exponent $4/3$ mentioned above can be replaced by $3/2$.

In this paper we study the bipartite case for general dimension, extending the results of 
\cite{Char,PdZ,Razdist}. We have the following main result.
\begin{theorem}   \label{thm:main}
For \(d \geq 2\), let \(\gamma_1,\gamma_2\) be a pair of constant-degree irreducible algebraic curves in \(\mathbb{R}^d\). Assume that $\gamma_i$ is neither contained in a hyperplane nor in a quadric surface in $\mathbb{R}^d$, for each $i=1,2$. Then, for every pair \(P_1 \subset \gamma_1,\; P_2 \subset \gamma_2\) of n-point sets, \(P_1 \times P_2\) spans \(\Omega(n^{3/2})\) distinct distances. 
\end{theorem}

Some special cases of the problem, where the curves are of one of the forms that are excluded in Theorem~\ref{thm:main}, are studied in Aldape et al.~\cite{ALPSV}.

\paragraph{Sketch of proof.} For the proof, we use tools from rigidity theory. Concretely, we use a theorem of Raz, Sharir, and De Zeeuw combined with the recent improvement of the exponent from Solymosi and Zahl (see Theorem~\ref{thm:2cases}), to show that in case $P_1\times P_2$ induce $o(n^{3/2})$ distinct distances, then the curves $\gamma_1$ and $\gamma_2$ support arbitrarily many equivalent realizations of the bipartite complete graph $K_{m,n}$, for $m,n$ arbitrarily large. Using the fact that $K_{m,n}$ is (generically) rigid in $\R^d$, and in view of the characterization of Bolker and Roth~\cite{BolkerRoth1980} for the non-generic realizations of $K_{m,n}$ (see Theorem~\ref{inf_rigid_condition}), we conclude that for each $i=1,2$ there exist arbitrarily many rigid motions of $\R^d$ that map the curve $\gamma_i$ to itself. We then show that such curves have derivatives with constant norm, and we apply D'Angelo and Tyson~\cite{helixTheorem} (see Corollary~\ref{cor:quadric} below) to show that $\gamma_1,\gamma_2$ must be contained in a quadric surface in this case, contradicting our assumptions. 

\paragraph{Organization of the paper.} In Section~\ref{sec:pre} we introduce the tools needed for the proof of our main Theorem~\ref{thm:main}. The proof of Theorem~\ref{thm:main} is then given in Section~\ref{sec:proofmain}. Some discussion and open problems are given in Section~\ref{sec:discussion}.

\section{Preliminaries and preparation for the proof}\label{sec:pre}

\subsection{Polynomials vanishing on Cartesian products}
We recall the following (special case of a) result of Raz, Sharir, and De Zeeuw~\cite{RazSharir16}, combined with the improved exponent from Solymosi and Zahl~\cite[Theorem 1.4]{SolyZahl}.

\begin{theorem} 
\label{thm:2cases} 
Let \( F \in \mathbb{R}[x,y,x]\) be an irreducible
polynomial of degree d, and assume that none of the derivatives \(\partial F / \partial x, \; \partial F/ \partial y,  \; \partial F / \partial z\) is identically zero. Then one of the following holds.
\begin{enumerate}[label=(\alph*)]
\item \label{case1} For all \(A,B,C \subset \mathbb{R}\; with\; |A| = |B| = n, \;|C| =\Omega(n)\), we have
\[|Z(F)\cap (A\times B\times C)| = O(|C|^{4/7}n^{8/7}),\]
where the constant of proportionality depends only on d. 

\item \label{case2} There exists a one-dimensional subvariety \(Z_0 \subset Z(F)\), such that for every \(v \in Z(F)\setminus Z_0\), there exist open intervals \(I_1, I_2, I_3 \subset \mathbb{R}\), and one-to-one real-analytic functions \(\varphi_i : I_i \to \mathbb{R} \), for \(i = 1,2,3\) with real-analytic inverses, such that \(v\in I_1 \times I_2 \times I_3\), and, for every \((x,y,z) \in I_1 \times I_2 \times I_3\):
\[(x,y,z) \in Z(F) \text{ if and only if } \varphi_3(z) = \varphi_1(x) -\varphi_2(y)\]
\end{enumerate}
\end{theorem}

\subsection{Graph rigidity}

We review the standard definitions from rigidity theory; for more details see e.g. Asimow and Roth~\cite{AsimowRoth1}. Let \(G: = (V,E)\) be a graph. A \emph{realization} \(\p\) of $G$ in \( \mathbb{R}^d\) is an embedding of the vertex set \(V = \{1,\ldots,|V|\}\) of $G$ in $\R^d$. That is, \[\p = (p_1,..,p_{|V|}) \in \mathbb{R}^d \times ... \times \mathbb{R}^d \cong \mathbb{R}^{d|V|}.\] 
A pair $(G,\p)$ of a graph $G$ and a realization $\p$ is called a \emph{framework}.

Define the \emph{edge function} of $G$
to be the map \(f_G: \mathbb{R}^{d|V|} \to \mathbb{R}^{|E|}\) given by \[(p_1,\ldots,p_{|V|})\mapsto \left(||p_i - p_j||^2\right)_{\{i,j\}\in E},\]
where we fixed some (arbitrary) order on the edges of $G$.
Note that \(f_G\) is a polynomial map.

For a given graph $G$, let $\p$ and $\q$ be a pair of realizations of $G$ in $\R^d$. We say that the corresponding frameworks, $(G,\p)$ and $(G,\q)$, are \emph{equivalent} if \(f_{G}(\p) = f_{G}(\q)\). Namely, if $\|p_i-p_j\|=\|q_i-q_j\|$ for every edge $\{i,j\}\in E$. We say that the frameworks $(G,\p)$ and $(G,\q)$ are \emph{congruent} if $\|p_i-p_j\|=\|q_i-q_j\|$ for every size-2 subset $\{i,j\}\subset V$ (here $\{i,j\}$ is not necessarily an edge in $G$). Equivalently, $(G,\p)$ and $(G,\q)$ are congruent if there exists a rigid motion $R$ of $\R^d$ such that $R(p_i)=q_i$, for every $i\in V$.

We say that \((G,p)\) is a \emph{rigid framework} if there exists a neighborhood \(U\) of \(\p\), such that for every \(\q\in U\), if $(G,\p)$ and $(G,\q)$ are equivalent, then they are necessarily also congruent. Equivalently, there exists a neighborhood $U$ of $\p$ such that
$$
f_G^{-1}(f_G(\p))\cap U=f_{K_{|V|}}^{-1}(f_{K_{|V|}}(\p))\cap U,$$
where $K_{|V|}$ is the complete graph on $|V|$ vertices.

Let \(\vartheta(d)\) be the \(d(d+1)/2\)-dimensional manifold of isometries of \(\mathbb{R}^d\). Define \(F_{\p}: \vartheta(d) \to \mathbb{R}^{d|V|}\) by 
$$S\mapsto (S(p_1), .. , S(p_{|V|})).$$ 
Observe that, for a realization \(\p = (p_1,...,p_{|V|}) \in \mathbb{R}^{d|V|}\) of $G$, the framework \((G,\p)\) is a rigid if and only if there exists a neighborhood \(U\) of \(\p\), such that 
$$
f^{-1}_{G}(f_G(\p)) \cap U = \im F_{\p} \cap U. 
$$
Assume in addition that the points of $\p$ affinely-span $\R^d$. 
Note that in this case $\im F_\p$ is a $d(d+1)/2$-manifold and that $\p\in \im F_\p$. Let \(T_\p\) denote the tangent space of the manifold $\im F_\p$ at the point \(\p\).

Let $df_G(\p)$ denote the Jacobian matrix of the function $f_G$, viewed as a linear transformation $df_G(\p): \mathbb{R}^{d|V|} \to \mathbb{R}^{|E|}$. The following lemma is from Roth~\cite{Roth81}; for completeness we provide the proof. 

\begin{lemma}\label{Tpsubsetker}
Let  \((G,\p)\) be a framework of a graph $G=(V,E)$. Then $T_\p \subseteq \ker df_G(\p)$.
\end{lemma}
\begin{proof}[Proof (based on Roth~\cite{Roth81}).]
Consider the matrix $df_G(\p)$. We may associate each of its rows with an edge of $G$, and every $d$ consecutive columns with a vertex of $G$. One can then verify that the $(\{i,j\},k)$ entry (corresponding to a $1\times d$ submatrix of $d f_G(\p)$) is given by
$$
\big(df_G(\p)\big)_{\{i,j\}, k} = 
\begin{cases}
2(p_i-p_j) & \text{if } k = i,\; \{i,j\}\in E\\
2(p_j-p_i) & \text{if } k = j,\; \{i,j\}\in E\\
0 & otherwise
\end{cases}
$$
Observe that for any $\mu=(\mu_1,\ldots,\mu_{|V|}) \in \mathbb{R}^{d|V|}$ we have that 
$d(f_G(\p))\mu = 0$
if and only if 
$$
(p_i-p_j)\cdot\mu_i + (p_j-p_i)\cdot\mu_j  = 0,
$$
or 
\begin{equation}
(p_i-p_j)\cdot(\mu_i-\mu_j) = 0. \label{eq:kernel}
\end{equation}
for every $\{i,j\}\in E$.

Let $\mu \in T_\p$. By definition, there exists a smooth path $x = (x_1, ..., x_{|V|}) : [0,1] \to \im F_\p$ such that $x(0) = \p$ and $\mu = x'(0)$. 
Then, by the definition of $F_\p$, we have that $x(t)$ is that image of $\p$ under some rigid motion of ${\mathbb R}^d$, for every $t\in[0,1]$. That is, for every $\{i,j\} \in \binom{V}{2}$ and every $t \in [0,1]$, we have
$$
||{x_i}(t) - {x_j}(t)||^2  = ||p_i- p_j||^2.
$$
Differentiating both sides with respect to $t$, and evaluating at $t = 0$, we get 
$$
({x_i}(0) - {x_j}(0))\cdot({x_i}'(0) - {x_j}'(0)) = ( {p_i} - {p_j})\cdot(\mu_i - \mu_j) = 0, 
$$
for every $\{i,j\}\in \binom{V}{2}$.
In particular, $\mu$ satisfies \eqref{eq:kernel}. Hence $T_\p \subseteq \ker df_G(\p),$ as needed. 
\end{proof}

\begin{defin}
We call a vector in \(\ker df_G(\p)\) an \emph{infinitesimal motion} of \((G,\p)\).
\end{defin}

\begin{defin}
Let \(\p\) be a realization of \(G = (V,E)\) in \(\mathbb{R}^d\). We say that the framework \((G,\p)\) is \emph{infinitesimally rigid} if its infinitesimal motions arise from rigid motions, that is, \( \ker df_G(\p)=T_\p\).
If the vertices of \((G,\p)\) do not lie on a hyperplane of  \(\mathbb{R}^d\), this condition is equivalent to \( \dim \ker df_G(\p) = d(d+1)/2\)
\end{defin}

Let \(X, Y\) be smooth manifolds and \(f: X \to Y\)
a smooth map.
Let \(k = \max{\{\rank df(x) : x \in X\}}\). 
We say that a point \(x \in X\) is a \emph{regular point} of \(f\) if \(\rank df(x) = k\), and a singular point otherwise.

\begin{theorem}[Asimow and Roth~\cite{AsimowRoth1}, Roth~\cite{Roth81}] \label{thm:rigidIffInfrigid}
Let \(G = (V,E)\) be a graph and $\p$ a realization of $G$. Suppose that $\p$ is a regular point in $\R^{d|V|}$ and that the points of $\p$ do not lie on a common hyperplane in $\R^d$.
Then
\((G,\p)\) is rigid if and only if \((G,\p)\) is infinitesimally rigid. 
\end{theorem}
\begin{proof}
Let $G$ and $\p$ be as in the statement.
Let 
$k = \max\{\rank df_G(x) : x\in \mathbb{R}^{d|V|}\}.$ 
As $\p$ is a regular point, we have that $\rank df_G(\p) = k$; moreover,  we have $\rank df_G(\q) = k$, for every $\q$ in some small neighborhood $U$ of $\p$. By the implicit function theorem (see e.g. \cite[p. 32]{Book}), assuming $U$ is sufficiently small, we have that $f_G^{-1}(f_G(\p))\cap U$ is a $(d|V|-k)$-dimensional smooth manifold. 

By the definition of rigidity, we have that \((G,\p)\) is rigid if and only if there is a neighborhood \(W\) of \(\p\) such that 
$$
\im F_\p \cap W = f^{-1}_G(f_G(\p)) \cap W.
$$ 
Note that $\im F_\p \subset  f^{-1}_G(f_G(\p))$ because preserving distances between all pairs of vertices, in particular preserves edge lengths. 
This implies that $(G,\p)$ is rigid if and only if $\dim \im F_\p = \dim f_G^{-1}(f_G(\p))$ near $\p$, that is, if and only if 
\begin{equation}\label{dimkrigid}
d(d+1)/2=d|V|- k.
\end{equation}

Then $(G,\p)$ is rigid if and only if 
\begin{align}\label{eq:dimenssionUsage}
k=\rank df_G(\p) 
&= d|V| - \dim \ker df_G(\p)\nonumber\\
&\leq d|V| - \dim T_\p \\
&= d|V|  -  d(d+1)/2\nonumber\\
&=k,\nonumber
\end{align}
    where for the second line we use Lemma~\ref{Tpsubsetker}, for the third line we use the fact that $\dim T_\p=d(d+1)/2$ and the last line is by \eqref{dimkrigid}.

We conclude that $(G,\p)$ is rigid if and only if the inequality \eqref{eq:dimenssionUsage} is in fact equality, that is, if and only if 
$\dim \ker df_G(\p)  = \dim T_\p$, which completes the proof. \end{proof}

\begin{theorem}[Asimow and Roth, Section 3, Theorem~\cite{AsimowRoth2}]\label{thm:notRigirToNotInfinitesimalRigid}
Let $G = (V,E)$ be a graph, \(\p = (p_1,...,p_{|V|}) \in \mathbb{R}^{d|V|}\) a realization in which the vertices do not lie on a common hyperplane in \(\mathbb{R}^d\), and \((G,\p)\) the associated framework.
If \((G,\p)\) is infinitesimally rigid then \((G,\p)\) is rigid and \(\p\) is regular.
\end{theorem}
\begin{proof}
Let $G$ and $\p$ be as in the statement, and assume that the framework $(G,\p)$ is infinitesimally rigid. By definition, this means that \(T_\p = \ker df_G(\p)\) and, since by assumption the points of $\p$ do not lie on a common hyperplane, we have in particular that 
\begin{align}
\rank df_G(\p) 
&= d|V| - \dim \ker df_G(\p)\nonumber\\
&=d|V| - d(d+1)/2.  \label{rankdfgp}  
\end{align}

Let \(q \in \mathbb{R}^{d|V|}\) be a regular point of $f_G$, and assume that the points of $\q$ do not lie on a common hyperplane in \(\mathbb{R}^d\).

So 
\begin{align*}
\rank df_G(\q) 
&\geq \rank df_G(\p)\\
&= d|V| - d(d+1)/2 \\ 
&=d|V|-\dim T_\q\\
&\geq d|V|-\dim\ker df_G(\q)\\
&=\rank df_G(\q),
\end{align*}
where for the second line we used \eqref{rankdfgp}.
We conclude that all inequalities above are in fact equalities, and in particular $\rank df_G(\p)=\rank df_G(\q)$. 
Hence $\p$ is a regular point.
By Theorem~\ref{thm:rigidIffInfrigid} we get also that $(G,p)$ is rigid.
\end{proof}

Note that the opposite direction of this theorem is true due to Theorem~\ref{thm:rigidIffInfrigid}. 

\subsection{Rigidity of complete bipartite graphs}
We review some properties of complete bipartite graph from Bolker and Roth~\cite{BolkerRoth1980}.

Let $A$ be a set of $k$ points in \(\mathbb{R}^d\). Write \(A = (a_1, a_2,..,a_k)\in \mathbb{R}^d \times \mathbb{R}^d \times ... \times \mathbb{R}^d \in \mathbb{R}^{kd}\), for some arbitrary order. 
The \emph{affine span} of $A$, denoted by $\Bar{A}$ is  the set 
$$
\Bar{A} = \left\{\sum_{i = 1}^k {\mu_i}{a_i} ~:~ \sum_{i = 1}^k {\mu_i} = 1\right\}
$$
Let 
$$
D(A) = 
\left\{(\lambda_1, \ldots, \lambda_k) \in \mathbb{R}^k  ~:~\sum_{i = 1}^k {\lambda_i}{a_i} = 0,\; \sum_{i = 1}^k{\lambda_i} = 0\right\}
$$
be the space of \emph{affine dependencies} of $A$.

By the dimension theorem for vector spaces, we have
\begin{equation} \label{calc: dimD(A)}
    \dim D(A) +  \dim \Bar{A} = k-1.
\end{equation}

\begin{defin}[\cite{BolkerRoth1980}]
Let $G=(V,E)$ be a graph and let $\p$ be a realization of $G$.
A \emph{stress} of the framework $(G,\p)$ is a function \(\omega : E \to \mathbb{R}\) such that 
\begin{equation*}
    \sum_{\left\{j : \{i,j\}\in E\right\}}{\omega_{\{i,j\}}(p_i-p_j)} = 0,
\end{equation*}
for each $i\in V$.
\end{defin}

Consider the complete bipartite graph  $G  = K_{m,n}$, and let 
$$
R = \left\{(A, B) = (a_1, ..., a_m, b_1, \ldots, b_n) : a_i,b_j \in \mathbb{R}^d\right\}
$$
denote the manifold of all realizations of $K_{m,n}$ in \(\mathbb{R}^d\). Let $K(A,B)$ denote the framework associated with a realization $(A,B)\in R$.

Let \((A,B)\in R\),
then a stress of \(K(A,B)\) has the form of an \(m \times n\) matrix \(\omega := (\omega_{i,j})\) such that:
\begin{gather}
    \nonumber \sum_{j=1}^n{\omega_{i,j}(a_i-b_j)} = 0, \quad\forall 1\leq i \leq m,\\
    \nonumber \sum_{i=1}^m{\omega_{i,j}(b_j-a_i)} = 0, \quad\forall 1\leq j \leq n.
\end{gather}
Define \(\Omega:= \Omega(A,B)\) to be the vector space of stresses of a framework \(K(A,B)\).

We observe that stresses of $K(A,B)$ are linear dependencies among the $nm$ rows of $df_{K_{m,n}}(A,B)$, thus 
$$nm  - \dim \Omega (A,B) = \rank f_{K_{m,n}}(A,B) = (m+n)d-\dim \ker df_{K_{m,n}}(A,B).$$
Thus
\begin{equation}\label{calcDimKer}
    \dim \ker df_{K_{m,n}}(A,B) = \dim \Omega (A,B) + (m+n)d -mn.
\end{equation}

Recall that a \emph{quadratic polynomial} is a
polynomial of the form
\[q(\xi_1, ..., \xi_d) = \sum_{1\leq i \leq j \leq d}{\sigma_{ij}\xi_i\xi_j} + \sum_{1\leq i\leq d} {\sigma_i}{\xi_i}+ \sigma.\]
Note that $q$ is determined by its \(d(d-1)/2 + 2d +1 = d(d+3)/2 +1\) coefficients. The set of quadric polynomials in \(\mathbb{R}^d\) forms a vector space denoted by $Q_d$.
A \emph{quadric surface} is the set of zeros of a nontrivial (i.e., not identically zero) quadric polynomial.

\begin{defin}
Let $C = (c_1, ..., c_k)\in \mathbb{R}^{kd}$ be a $k$-tuple of  vectors in \(\mathbb{R}^d\), and assume that $\Bar{C} = \R^d$. Let 
$$
Q(C) = \{q\in Q_d ~:~ q(c) = 0~\text{for every $c\in C$} \}.
$$
\end{defin}

\begin{theorem}[{Bolker and Roth~\cite[Thm. 10]{BolkerRoth1980}}] \label{calc:dimOmega}
Let $(A,B)\in R$ be a realization of $K_{m,n}$ and put $C = (A \cap \Bar{B}, B \cap \Bar{A}) = (c_1, ...c_k) \in \mathbb{R}^{kd}$. Assume that $\bar{C}=\R^d$. Then
\begin{equation}
      \dim \Omega(A, B) = \dim D(A) \dim D(B) + \dim Q(C) + k - d( d + 3)/2 - 1.
\end{equation}
\end{theorem}

We conclude the following fact, which is a direct consequence of the analysis in Bolker and Roth~\cite{BolkerRoth1980}.
\begin{theorem}\label{inf_rigid_condition} Consider the complete bipartite graph $K_{m,n}$, and let $(A,B)\in R$ be a realization of $K_{m,n}$ in $\R^d$. Assume that $\Bar{A} = \Bar{B} = \mathbb{R}^d$. Then $K(A,B)$ is infinitesimally rigid unless its $m+n$ vertices lie on a quadric surface.
\end{theorem}
\begin{proof}
Let $(A,B)\in R$ be a realization of \(K_{m,n}\) in \(\mathbb{R}^d\) for which \(\Bar{A} = \Bar{B} = \mathbb{R}^d\). So \(C = (A \cap \Bar{B}, B \cap \Bar{A}) = (A, B) \). By \eqref{calc: dimD(A)}, we have 
\begin{align*}
\dim D(A) &= m-d-1\\ 
\dim D(B) &= n-d-1.
\end{align*}
By Theorem~\ref{calc:dimOmega}, we have 
\[
\dim \Omega(A, B) = (m-d-1)(n-d-1) + \dim Q(C) + (m+n) - d(d+3)/2 - 1.
\] 
Using \eqref{calcDimKer}, we get
\begin{align*}
\dim \ker df_{K_{m,n}} 
&= (m-d-1)(n-d-1) + \dim Q(C)+ m+n-d(d+3)/2 - 1+(m+n)d-mn\\
&= d(d +1)/2 +\dim Q(C)\\
&=\dim T_{(A,B)}+ \dim Q(C).
\end{align*}
So \(\dim T_{(A,B)} = \dim \ker df_{k_{m,n}}(A,B)\) if and only if \(\dim Q(C) = 0\). It thus follows that the $K(A,B)$ is infinitesimally rigid unless its \(m+n\) vertices lie on a quadric surface, as asserted.
\end{proof}

\subsection{Characterization of curves with constant-norm derivatives}\label{sec:constderiv}
A \emph{real analytic parameterized curve} in \(\mathbb{R}^d\) is a real analytic map \(\gamma\) from an interval \(I \subseteq \mathbb{R}\) to $\mathbb{R}^d$. We often identify \(\gamma\) with its image in $ \mathbb{R}^d $.
For \(j\geq 0\), we denote by \(D^j\gamma\) the $j$th order
derivative of \(\gamma\) with respect to the parameter t (where \(D^0\gamma := \gamma\)).

\begin{defin}
For \(k \geq 0\), let
$Q_k = Q^{(d)}_k$ 
be the collection of real analytic parameterized curves \(\gamma \subset \mathbb{R}^d\) such that \(D^j\gamma\) has
constant Euclidean norm for each \(j \geq k.\)
\end{defin}

The following characterization is due to D’Angelo and Tyson~\cite{helixTheorem}. In the statement of
the theorem we use the standard notation
$\exp(B)=\sum_{k\ge 0}\frac{B_k}{k!}$, for a matrix $B$.

\begin{theorem}[{{\bf D'Angelo and Tyson~\cite[Theorem 3.7]{helixTheorem}}}]\label{Q0}
Let $\gamma:I\to\R^d$ be a real analytic parameterized curve, for some nonempty open 
interval $I\subset\R$. Assume that $\gamma\in Q_0$. Then, 
up to a rigid motion,
$\gamma$ is given by
\begin{equation}\label{Q0eq}
\gamma(t)=(\exp(At)v,w)\in\R^{2k}\times\R^{d-2k} ,\qquad \text{for $t\in I$} ,
\end{equation}
where $A$ is a real skew-symmetric $2k\times 2k$ invertible matrix, $v\in \R^{2k}$, $w\in \R^{d-2k}$ and $0\le k\le d/2$ is some integer.
\end{theorem}
Observe that for $\gamma\in Q_1$ one has $D^1\gamma\in Q_0$. Thus, applying Theorem~\ref{Q0}, 
and integrating the resulting representation (\ref{eqq0}), we have the following result.

\begin{corollary}[{D’Angelo and Tyson~\cite[Corollary~3.8]{helixTheorem}}]
\label{helixForm}
Let \(\gamma: I \to \mathbb{R}^d\)
be a real analytic
parameterized curve, for some nonempty open interval \(I \subset \mathbb{R}\). Assume that \(\gamma \in Q_1\). Then,
up to a rigid motion,
\(\gamma\) is given by:
\begin{equation*}
    \gamma(t) = (\exp(At)v+v_0,tw+w_0) \in \mathbb{R}^{2k}\times \mathbb{R}^{d-2k}, \quad \text{for $t\in I$}
\end{equation*}
where A is a real skew-symmetric $2k \times 2k$ invertible matrix, $v, v_0 \in \mathbb{R}^{2k}$, $w, w_0 \in \mathbb{R}^{d-2k}$ and $0 \leq k \leq d/2 $ is some integer.
\end{corollary}

Recall that the spectral decomposition of $2k\times 2k$ invertible real skew-symmetric matrix $B$ is $B=U\Sigma U^{-1}$ (see e.g. ~\cite[p. 293]{skew_symetric_decomp}), where  $U$ is an orthogonal linear transformation
and $\Sigma$ is a block diagonal matrix of the form
$$
\Sigma = 
\begin{bmatrix}
\begin{matrix}0 & \lambda_1\\ -\lambda_1 & 0\end{matrix} &  0 & \cdots & 0 \\
0 & \begin{matrix}0 & \lambda_2\\ -\lambda_2 & 0\end{matrix} &  & 0 \\
\vdots &  & \ddots & \vdots \\
0 & 0 & \cdots & \begin{matrix}0 & \lambda_k\\ -\lambda_k & 0\end{matrix} 
\end{bmatrix}
$$
for some real $\lambda_1,\ldots,\lambda_k$.
Some manipulations then show that $\exp(B) = U\exp(\Sigma) U^{-1}$, and
$$
\exp(\Sigma) = 
\begin{bmatrix}
\begin{matrix}\cos\lambda_1 & \sin\lambda_1\\ -\sin\lambda_1 & \cos\lambda_1\end{matrix} &  0 & \cdots & 0 \\
0 & \begin{matrix}\cos\lambda_2 & \sin\lambda_2\\ -\sin\lambda_2 & \cos\lambda_2\end{matrix} &  & 0 \\
\vdots &  & \ddots & \vdots \\
0 & 0 & \cdots & \begin{matrix}\cos\lambda_k & \sin\lambda_k\\ -\sin\lambda_k & \cos\lambda_k\end{matrix} 
\end{bmatrix} .
$$

We conclude that, for a curve $\gamma \in Q_1$, Corollary~\ref{helixForm} implies that, up to a rigid motion,
we have
\begin{equation} \label{eqq0}
\gamma(t)=(\rho_1\cos(\lambda_1t +\theta_1),\rho_1\sin(\lambda_1t+\theta_1),\ldots,
\rho_k\cos(\lambda_kt+\theta_k),\rho_k\sin(\lambda_kt+\theta_k),tw),
\end{equation} 
for some $\rho_1,\ldots,\rho_k \in\R\setminus\{0\}$, $\lambda_1,\ldots,\lambda_k\in\R\setminus \{0\}$, and $\theta_1,\ldots,\theta_k\in [0,2\pi)$.

Finally, for the proof of our main Theorem~\ref{thm:main}, we will need the following simple corollary.
\begin{corollary}\label{cor:quadric}
Let \(\gamma: I \to \mathbb{R}^d\)
be a real analytic
parameterized curve, for some nonempty open interval \(I \subset \mathbb{R}\). Assume that \(\gamma \in Q_1\). Then $\gamma$ is contained in a quadric surface in $\R^d$.
\end{corollary}
\begin{proof}
As argued above, up to a rigid motion, such $\gamma$ has a parameterization of the form \eqref{eqq0}. 
In particular, if $k=0$ then $\gamma$ is contained in a line, and otherwise, denoting the coordinates in $\R^d$ by $(x_1,\ldots,x_d)$, every point on $\gamma$ satisfies the equation
$x_1^2+x_2^2=\rho_1^2$. Thus $\gamma$ is contained in a quadric surface, as asserted.
\end{proof}

\begin{remark}
Note that, in Section~\ref{sec:proofmain}, we apply Corollary~\ref{cor:quadric} to a curve $\gamma$ that is an arc of an irreducible \emph{algebraic} curve. In this case, one can say more about the structure of $\gamma\in Q_1$. Concretely, that it is either a line or in fact in $Q_0$. See Charalambides~\cite[Definition 1.5 and Lemma 7.4]{Char}.
\end{remark}

\section{\textbf{Proof of Theorem \ref{thm:main}}}\label{sec:proofmain}
For each \(k = 1,2\), we may assume without loss of generality that the curve \(\gamma_k \subset \mathbb{R}^d\) is given by the system
\[g_i^k(x_1,...,x_d) = 0,\qquad i = 1,..,d-1,\]
where each of \(g_i^k\) is an irreducible constant-degree d-variate real polynomial. For every pair of points \(\boldsymbol{x} = (x_1,..,x_d)\in \gamma_1\) and \(\boldsymbol{y} = (y_1,..,y_d)\in \gamma_2\) at distance \(r\geq0\), we have:
\[g^1_i(\boldsymbol{x}) = 0, \qquad i = 1,...,d-1,\]
\[g^2_i(\boldsymbol{y}) = 0, \qquad i = 1,...,d-1,\]
\[||\boldsymbol{x}-\boldsymbol{y}||^2-r^2 = (x_1-y_1)^2+(x_2-y_2)^2+...+(x_d-y_d)^2-r^2 = 0\]

The system above defines an algebraic variety \(V\) in \(\mathbb{R}^{2d+1}\) with coordinates \((\boldsymbol{x}, \boldsymbol{y}, r)\), of degree that depends on \(\deg\gamma_1\), \(\deg\gamma_2\) and \(d\). Note that \(V\) is two-dimensional. Indeed,  
the projection of \(V\) to the first \(2d\) coordinates of \(\mathbb{R}^{2d+1}\) is exactly \(\gamma_1\times\gamma_2\), which is two-dimensional, and each point in this projection, has a finite pre-image in \(V\).

By applying a generic isometry of \(\mathbb{R}^d\), 
we may assume without loss of generality that
each of $\gamma_1,\gamma_2$ has finite fibers with respect to the projection of $\R^d$ onto its first coordinate. Let \(\pi : \mathbb{R}^{2d+1}\to\mathbb{R}^3\) denote the projection of $\R^{2d+1}$ to the coordinates \((x_1, y_1, r)\). Then  \(W:=\pi(V)\) is a semi-algebraic set and the pre-image in $V$ of a generic element of $W$ is finite. Moreover, since taking projection neither increases the degree of a variety nor its dimension, the Zariski closure \(Cl(W)\) of \(W\) is two-dimensional
and its degree is a constant that depends only on  \(\deg\gamma_1\), \(\deg\gamma_2\) and \(d\).

Then there exists a trivariate constant-degree real polynomial \(F(x, y, z)\) such that \(Z(F) =
Cl(W)\). We are now ready to apply Theorem 1. By construction, 
it is easy to see that \(F\) depends non-trivially on each of its variables. In
what follows we assume that \(F\) is irreducible. 
Otherwise we apply the forthcoming analysis
to each irreducible component of \(F\) separately. Hence, \(F\) satisfies one of the properties (a)
or (b) of Theorem \ref{thm:2cases}.

\paragraph{Case 1: $F$ satisfies property (a) of Theorem~\ref{thm:2cases}.}
Let $P_1\subset \gamma_1$ and $P_2\subset \gamma_2$ be finite sets, each of size $n$.
Recall our assumption that
each of $\gamma_1,\gamma_2$ has finite fibers with respect to the projection of $\R^d$ onto its first coordinate.
Let $A,B$ denote the set of first coordinates of the points in $P_1,P_2$, respectively.
Then each of $A$, $B$ is of size $\Theta(n)$.
Let $C$ denote the set of distances spanned by $P_1\times P_2$. 
By property (a) in Theorem~\ref{thm:2cases}, we have 
\begin{equation}\label{upper}
|Z(F)\cap (A\times B\times C)|=O(|C|^{4/7}n^{8/7}).
\end{equation}

Next we obtain a lower bound on $|Z(F)\cap (A\times B\times C)|$. 
For this, note that for every pair $(a,b)\in A\times B$ there exists $c\in C$ such that $F(a,b,c)=0$.
Indeed, let $p\in P_1$ be a point whose first coordinate is $a$, and let $q\in P_2$ be 
a point whose first coordinate is $b$, and let $c:=\|p-q\|$. By definition, $c\in C$ and $(p,q,c)\in V$.
Hence, $(a,b,c)=\pi(p,q,c)\in W\subset Z(F)$, which means in particular that $F(a,b,c)=0$. Thus, 
\begin{equation}\label{lower}
|Z(F)\cap (A\times B\times C)|=\Omega(|A||B|).
\end{equation}
(If $F$ is reducible, this must hold for at least one irreducible component of $F$.)
Combining \eqref{upper} and \eqref{lower}, we conclude that 
$$
|C|=\Omega(n^{3/2}),
$$
which completes the proof of Theorem~\ref{thm:main} for this case.

\paragraph{Case 2: F satisfies property (b) of Theorem~\ref{thm:2cases}.}  We will prove that in this case each of $\gamma_1$ and $\gamma_2$ is contained in a quadric surface. 

Let \(u = (x_1, y_1, z_1) \in (Z(F)\setminus Z_0) \cap W\),
where \(Z_0\) is the exceptional subset of \(Z(F)\) provided in property (b). Since \(Z(F)={\rm Cl}(W)\) and $Z_0$ is lower-dimensional, we have that the intersection \((Z(F)\setminus Z_0) \cap W = W \setminus Z_0\) is nonempty, and hence such an element \(u\) exists. Let \(U = I_1\times I_2 \times I_3\) be the open neighborhood of \(u\) provided in Theorem~\ref{thm:2cases} property \ref{case2}, such that
\[Z(F) \cap U  = \{(x,y,z)\in U | \varphi_3(z) = \varphi_1(x)-\varphi_2(y)\}\] 
for suitable real-analytic invertible functions \(\varphi_i : I_i \to \mathbb{R}\).

By shrinking \(U\), if needed, we may assume that \(Z(F)\cap U \subset W\). That is, we are in the following situation:
For each \(i = 1, 2\), we have some local parameterizations \(\alpha_i :I_i \to \gamma_i \subset \mathbb{R}^d\) of an arc of \(\gamma_i\), with the property that
\[||\alpha_1(x) - \alpha_2(y)|| = \varphi_3^{-1}(\varphi_1(x) - \varphi_2(y)).\]
By a proper change of variables, we may assume, without loss of generality, that in fact
\begin{equation}\label{distanceh}
||\alpha_1(x) - \alpha_2(y)|| = h(x-y)\,
\end{equation}
for some invertible analytic function \(h\) defined over the image of \((x,y)\to x-y\),  for \(x,\, y\) in
some (common) open interval \(I \subseteq \mathbb{R}.\)

\begin{lemma}\label{lemma:infCopies} Let \(m, n \geq d+1\). Then there exists $\delta_0>0$ such that the following holds. With each $0\le \delta\le \delta_0$, we may associate a realization $(A_\delta, B_\delta)= (a_{\delta,1}, ...a_{\delta,m}, b_{\delta,1},..., b_{\delta,n})$ of the complete bipartite graph $K_{m,n}$, such that 
the affine span of \((a_{\delta,1}, ..., a_{\delta, {d+1}})\) and of \((b_{\delta,1}, ..., b_{\delta,{d+1}})\) is \(d\)-dimensional and the frameworks $(K_{m,n}, (A_\delta, B_\delta))$ are all equivalent to each other.
\end{lemma}
\begin{proof}
By assumption, none of \(\gamma_1\) and \(\gamma_2\) is contained in a hyperplane in \(\mathbb{R}^d\), and so, using  B\'ezout's theorem, also none of \(\alpha_1\) and \(\alpha_2\) is contained in one. 
Thus, for each $i=1,2$,  there exist \(d+1\) points on \(\alpha_i\) such that their affine span is $d$-dimensional.

Let $A_0=\{\alpha_1(x_{0,1}), \ldots\alpha_1(x_{0,m})\}$ be $m$ points on $\alpha_1$, with $m\ge d+1$, such that \(\{x_{0,1}\leq x_{0,2} \leq \ldots \leq x_{0,m}\}\subset I\) and the affine span of $\alpha_1(x_{0,1}), \ldots\alpha_1(x_{0,d+1})$ is $\mathbb{R}^d$. Similarly, let $B_0=\{\alpha_2(y_{0,1}), \ldots\alpha_2(y_{0,n})\}$ be $n$ points on $\alpha_2$, for $n\ge d+1$, such that \(\{y_{0,1}\leq y_{0,2} \leq \ldots \leq y_{0,n}\}\subset I\) and the affine span of $\alpha_2(y_{0,1}), \ldots\alpha_2(x_{0,d+1})$ is $\mathbb{R}^d$.
We consider the realization of $K_{m,n}$ given by $(A_0, B_0)$.

Given $\delta_0>0$, define for every $0\leq \delta \leq \delta_0$
\begin{align*}
  x_{\delta, i} &= x_{0,i}+ \delta, \quad i=1,\ldots m\\
  y_{\delta, j} &= y_{0,j}+ \delta, \quad j=1,\ldots, n  
\end{align*}
and respectively
\begin{align*}
  a_{\delta, i} &= \alpha_1(x_{\delta, i}),  \quad i=1,\ldots m\\
  b_{\delta, j} &= \alpha_2(y_{\delta, j}),  \quad j=1,\ldots n.
\end{align*}
We may choose $\delta_0>0$ sufficiently small so that $x_{\delta, i}, y_{\delta,j}\in I$ for every $i\in\{1,\ldots, m\}$, $j\in\{1,\ldots,n\}$, $\delta\in[0,\delta_0]$, and so that each of the sets $\{a_1,\ldots, a_{d+1}\}$  and $\{b_1,\ldots,b_{d+1}\}$ affinely-spans $\mathbb{R}^d$.

For each $\delta\in [0,\delta_0]$, let $A_\delta:=\{a_{\delta, 1}, \ldots a_{\delta,m}\}$, $B_\delta:=\{b_{\delta, 1}, \ldots b_{\delta,n}\}$ and consider the realization $(A_\delta,B_\delta)$ of $K_{m,n}$ in $\R^d$.

We claim that the frameworks $(K_{m,n}, (A_\delta,B_\delta))$, for $0\le \delta\le \delta_0$, are equivalent. 
Indeed, the coordinate of the vector $f_{K_{m,n}}(A_\delta,B_\delta)$ corresponding to the edge $\{i,j\}$ in $K_{m,n}$ is
\begin{align*}
\left(f_{K_{m,n}}(A_\delta,B_\delta)\right)_{\{i,j\}} &= ||a_{\delta, i} - b_{\delta, j}||^2 \\
&=||\alpha_1(x_{0,i}+ \delta)- \alpha_2(y_{0,j}+ \delta)||^2 \\
&=(h(x_{0,i} - y_{0,j}))^2,
\end{align*}
where for the last equality we used \eqref{distanceh}.
Since the length is independent of $\delta$, this proves our claim and therefore completes the proof of the lemma.
\end{proof}

Let $(A_\delta,B_\delta)$, for $0\le\delta\le \delta_0$ be the realizations given by Lemma~\ref{lemma:infCopies}. We observe that if, for some $\delta\in[0,\delta_0]$, the framework $(K_{m,n}, (A_\delta,B_\delta))$ is not rigid, then $\gamma_1\cup\gamma_2$ is contained in a quadric surface, contradicting our assumption. Indeed, in this case, $(K_{m,n}, (A_\delta,B_\delta))$ is in particular not infinitesimally rigid (see Theorem~\ref{thm:notRigirToNotInfinitesimalRigid}), and then by Corollary \ref{inf_rigid_condition} we have that the points in $A_\delta\cup B_\delta$ are contained in some quadric surface $Q$. Choosing $m,n>0$ sufficiently large (depending on the degrees of $\gamma_1,\gamma_2$), B\'ezout's theorem implies that both $\gamma_1,\gamma_2$ are contained in $Q$.

Therefore, in what follows, we may assume that $(K_{m,n}, (A_\delta,B_\delta))$ is rigid, for every $\delta\in[0,\delta_0]$. Moreover, as the following lemma asserts, we may assume without loss of generality that those frameworks are in fact congruent.

\begin{lemma} \label{lemma:seqOfCongruent}
Let $m,n\ge 0$. Let $(K_{m,n},(A_\delta,B_{\delta}))$, for $\delta\in[0,\delta_0]$, be a family of equivalent frameworks. Assume in addition that $(K_{m,n}, (A_\delta,B_\delta))$ is rigid in $\mathbb{R}^d$, for each $\delta$. Then there exists $0<\delta'< \delta_0$, such that the subfamily of frameworks $(K_{m,n}, (A_\delta,B_\delta))$, $0\le \delta\le\delta'$ are all congruent to each other.
\end{lemma}
\begin{proof}
By assumption $(K_{m,n}, (A_{0}, B_{0}))$ is a rigid framework. That is, by definition, there exists a neighborhood $U$ of $(A_{0}, B_{0})$ such that, if $(A,B)\in U$ and $(K_{m,n}, (A,B))$ is equivalent to $(K_{m,n},(A_{0}, B_{0}))$, then the two frameworks are also congruent. Let $0< \delta'\le \delta_0$ be sufficiently small so that $(A_\delta,B_\delta)\in U$, for every $0\le \delta\le \delta'$. Then $(K_{m,n}, (A_\delta,B_\delta))$ is congruent to $(K_{m,n}, (A_0,B_0))$, for every $0\le \delta\le \delta'$, as needed.
\end{proof}

We are now ready to prove that each of $\gamma_1$ and $\gamma_2$ is contained in a quadric surface.

By Lemma~\ref{lemma:seqOfCongruent}, there is \(0 < \delta' \leq \delta_0\), such that the realizations \(\{(A_\delta, B_\delta) | 0 \leq \delta \leq \delta'\}\) are all congruent. In particular, for every \( c\in [0, \delta'] \) there is an isometry \(T_c : \mathbb{R}^d \to \mathbb{R}^d\) such that:
\begin{align}\label{isometryintersection}
T_c(\alpha_1(x_{0,i})) &=\alpha_1(x_{c,i}) = \alpha_1(x_{0,i}+c),\quad\text{for each $i=1,\ldots,m$ and}\\
T_c(\alpha_2(y_{0,j})) &=\alpha_2(y_{c,j})=\alpha_2(y_{0,j}+c)\quad \text{for each $j=1,\ldots,n$}.\nonumber
\end{align}
Note that $T_c(\gamma_1)$ is an irreducible algebraic curve with the same degree as $\gamma_1$. By \eqref{isometryintersection}, we have that $| T_c(\gamma_1)\cap\gamma_1|\ge m$. If $m$ is sufficiently large (namely, $m>\deg^2 \gamma_1$), then  B\'ezout's theorem implies that $T_c(\gamma_1)=\gamma_1$ (here we used our assumption that $\gamma_1$ is irreducible). 
Similarly, we conclude also that $T_c(\gamma_2)=\gamma_2$, by choosing $n$ sufficiently large (namely, $n> \deg^2\gamma_2$).

Next, we prove that, for each \( c\in [0, \delta'] \) and \(t\) in a neighborhood close enough to 0, one has
\begin{equation}\label{s=t}
    T_c(\alpha_1(x_{0,1} + t)) = \alpha_1(x_{c,1}+t).
\end{equation}

We claim that a neighborhood, $U_c$, of $\alpha_1(x_{0,1})$ is mapped by $T_c$ to a neighborhood $T_c(U_c)$ of $\alpha_1(x_{c,1})$ (a neighborhood  in $\alpha_1$). Indeed, by the continuity of $T_c$ and of $\alpha_1$, points in $\mathbb{R}^d$ that are close to $\alpha_1(x_{0,1})$ are mapped to points in $\mathbb{R}^d$ that are close to  $T_c(\alpha_1(x_{0,1}))=\alpha_1(x_{c,1})$. In addition, as has already been argued, points on $\gamma_1$ are mapped to points on $\gamma_1$, so our claim follows. 

As $\alpha_1$ is a bijection, we may write $U_c=\{\alpha_1(x_{0,1}+t)\mid t\in J_c\}$, where $J_c$ is some open interval containing 0. Similarly, we may write $T_c(U_c)=\{\alpha_1(x_{c,1}+s)\mid s\in J_c'\}$, where $J_c'$ is an open interval containing 0.

Using the fact that $T_c$ is a bijection, we have
\begin{equation}\label{ttos}
T_c(\alpha_1(x_{0,1}+t))=\alpha_1(x_{c,1}+s(t)),\quad \text{for every $t\in J_c$},
\end{equation}
where $s:J_c\to J_c'$ is a bijection.

Recalling that $T_c$ is an isometry, we have, for every $t\in J_c$, that
\begin{align}\label{eqqq}
||\alpha_1(x_{0,1}+t)- \alpha_2(y_{0,1})|| &=||T_c(\alpha_1(x_{0,1}+t))-T_c(\alpha_2(y_{0,1}))||\\ &=||\alpha_1(x_{c,1}+s(t)) - \alpha_2(y_{0,1}+c)||,    \nonumber
\end{align}
where for the last line we used \eqref{isometryintersection} and \eqref{ttos}.

Applying \eqref{distanceh}, we get
\begin{align*}
||\alpha_1(x_{0,1}+t)- \alpha_2(y_{0,1})||  &= h(x_{0,1} - y_{0,1} + t)
\end{align*}
and
\begin{align*}
||\alpha_1(x_{c,1}+s(t)) - \alpha_2(y_{0,1}+c)|| &= h(x_{c,1} - y_{0,1}+ s(t) - c)\\
&= h(x_{0,1} + c - y_{0,1}+ s(t) - c) \\
&= h(x_{0,1} - y_{0,1}+ s(t))
\end{align*}
In view of \eqref{eqqq}, this gives
$$h(x_{0,1} - y_{0,1} + t)=h(x_{0,1} - y_{0,1}+ s(t)),$$
which in turn implies that $t=s(t)$, for every $t\in J_c$, as $h$ is a bijection. This proves \eqref{s=t}.

 Let \(I_{\delta'} = [x_{0,1}, x_{0,1} +\delta']\) and \(\tilde\alpha_{1} = \alpha_1|_{I_{\delta'}}\), the restriction of $\alpha_1$ to the interval $I_{\delta'}$. We claim that
\begin{equation}\label{tildealpha1inQ1}
\tilde\alpha_{1}\in Q_1
\end{equation}
(see Section~\ref{sec:constderiv} for the definition of $Q_1$).
To prove this, we show that, for every $u\in I_{\delta'}$ fixed, we have
\begin{equation}\label{derivnormconsttilde}
\|\tilde\alpha^{(k)}_{1}(u)\|=\|\tilde\alpha^{(k)}_{1}(x_{0,1})\|,
\end{equation}
for every positive integer $k$.
Note that, by the definition of $\tilde \alpha_1$ it is equivalent to prove that 
\begin{equation}\label{derivnormconst}
\|\alpha^{(k)}_{1}(u)\|=\|\alpha^{(k)}_{1}(x_{0,1})\|,
\end{equation}
for every $u\in I_{\delta'}$ fixed.
Note also that every $u\in I_{\delta'}$ is of the form $x_{0,1}+c$, for some $c\in [0,\delta']$. That is, we need to prove
\begin{equation}\label{derivnormconstc}
\|\alpha^{(k)}_{1}(x_{0,1}+c)\|=\|\alpha^{(k)}_{1}(x_{0,1})\|,
\end{equation}
for every $c\in [0,\delta']$ fixed.

We first prove \eqref{derivnormconstc} for $k=1$. 
By \eqref{s=t}, for every \(t\in J_c\), we have
\begin{equation}\label{s=ttilde} 
{\alpha}_1(x_{c,1} + t)=T_c({\alpha}_1(x_{0,1} + t)),
\end{equation}
where $J_c$ is some small interval containing 0.
Recall that $T_c$ is an affine transformation; more precisely, $T_c=A_c+\vec{b}$, where $A_c$ is an orthogonal $d\times d$ matrix and $\vec{b}$ is a vector in $\mathbb{R}^d$. Taking derivative with respect to $t$ from both sides, \eqref{s=ttilde} becomes
\begin{align*}
     \alpha_{1}'(x_{c,1} + t) &=  (T_c)'(\alpha_1(x_{0,1} + t))\cdot{\alpha}'_1(x_{0,1} + t)\\
     &=A_c\cdot{\alpha}'_1(x_{0,1} + t)
\end{align*}

Since $A_c$ is an isometry, we get
\begin{equation*}
||\alpha_1'(x_{c,1} +t)|| = ||A_c\alpha_1'(x_{0,1} + t)||\\
    = ||\alpha_1'(x_{0,1} + t)||.    
\end{equation*}
Setting $t=0$, we conclude that 
\begin{equation*}
||\alpha_1'(x_{c,1})|| = ||\alpha_1'(x_{0,1})||,    
\end{equation*}
for every $c\in [0,\delta']$. This proves \eqref{derivnormconstc} for $k=1$.

For $k\ge 2$, we note that
\begin{equation}\label{kderiv}
    \alpha_{1}^{(k)}(x_{c,1} + t)=A_c\cdot \alpha_{1}^{(k)}(x_{0,1} + t).
\end{equation}
Indeed, we prove by induction on $k$. Assume that 
$$\alpha_{1}^{(k-1)}(x_{c,1} + t)=A_c\cdot \alpha_{1}^{(k-1)}(x_{0,1} + t).$$ Then taking derivative from both sides with respect to $t$, we get
$$\alpha_{1}^{(k)}(x_{c,1} + t)=A_c\cdot \alpha_{1}^{(k)}(x_{0,1} + t).$$ 
Claim \eqref{derivnormconstc} now follows also for $k\ge2$, by combining \eqref{kderiv}, setting $t=0$, with the fact that $A_c$ is an isometry.

Thus we have proved \eqref{derivnormconstc} which implies \eqref{tildealpha1inQ1}.
By a symmetric argument, applied to $\tilde\alpha_2=\alpha_2\mid_{I_\delta'}$, one also has $\tilde\alpha_2\in Q_1$.

Applying Corollary~\ref{cor:quadric}, each of $\tilde\alpha_1$ and $\tilde \alpha_2$ is contained in a quadric surface. By B\'ezout's theorem this implies that the same holds for each of $\gamma_1,\gamma_2$. This completes the proof of the theorem. 
\hfill$\square$

\section{Discussion}\label{sec:discussion}
Our main theorem asserts the following: If $P_1\subset \gamma_1$ and $P_2\subset \gamma_2$ and the number of distinct distances spanned by $P_1\times P_2$ is ``small", then each of $\gamma_1$ and $\gamma_2$ must be contained in a quadric surface. One may wish to find a more concrete characterization of those ``special" $\gamma_1$ and $\gamma_2$ that allow a small number of distinct distances. In the plane Pach and De Zeeuw showed that $\gamma_1,\gamma_2$ must be either a pair of lines or a pair of concentric circles. In higher dimensions, for the non-partite case, the special curves are known to be algebraic helices (see Charalambides~\cite{Char} and also Raz~\cite{Razdist}). 

\begin{conj}\label{conjmain}
For \(d \geq 2\), let \(\gamma_1,\gamma_2\) be a pair of constant-degree irreducible algebraic curves in \(\mathbb{R}^d\).
Assume that neither of $\gamma_1$, $\gamma_2$ is contained in a hyperplane in $\mathbb{R}^d$.
Then, for every pair \(P_1 \subset \gamma_1,\; P_2 \subset \gamma_2\) of n-point sets, \(P_1 \times P_2\) spans \(\Omega(n^{3/2})\) distinct distances, unless each of $\gamma_1$ and $\gamma_2$ is an algebraic helix.
\end{conj}

As a step towards proving the above conjecture, we believe that the following characterization, which can be viewed as a certain extension of Bolker and Roth result (\cite{BolkerRoth1980}), should hold:
\begin{conj}\label{conjBR}
Consider the complete bipartite graph $K_{m,n}$, and let $(A,B)\in R$ be a realization of $K_{m,n}$ in $\R^d$. 
Then $K(A,B)$ is rigid  in $\R^d$, unless either the points in $A$ or the points in $B$ are contained in a hyperplane in $\R^d$.
\end{conj}
Note that in Conjecture~\ref{conjBR} we consider to the notion of \emph{rigidity}, unlike in Bolker and Roth~\cite{BolkerRoth1980} where the notion of \emph{infinitesimal rigidity} is considered (for infinitesimal rigidity, Bolker and Roth's results is an if and only if and cannot be strengthened). Note also that, by following the outline of the proof in Section~\ref{sec:proofmain}, Conjecture~\ref{conjBR} implies Conjecture~\ref{conjmain}.

\end{document}